\documentclass[12pt]{amsart}

\theoremstyle{plain}

\newtheorem{Thm}{ }[section]

\usepackage{amsfonts,ifthen,epsfig,amssymb}

\newcommand{\apic}[1]{\epsfig{figure=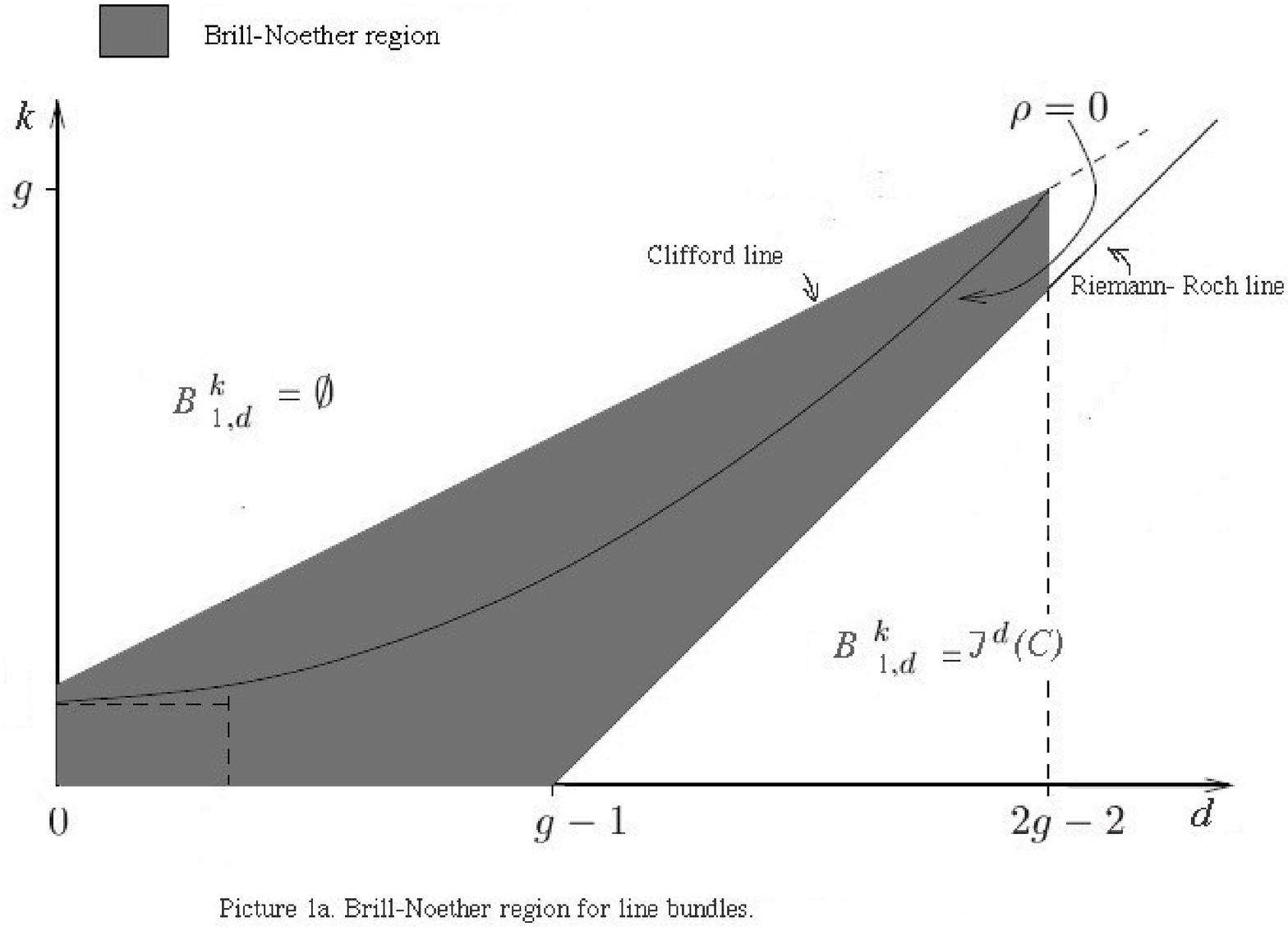,width=#1pt}}

\newcommand{\bpic}[1]{\epsfig{figure=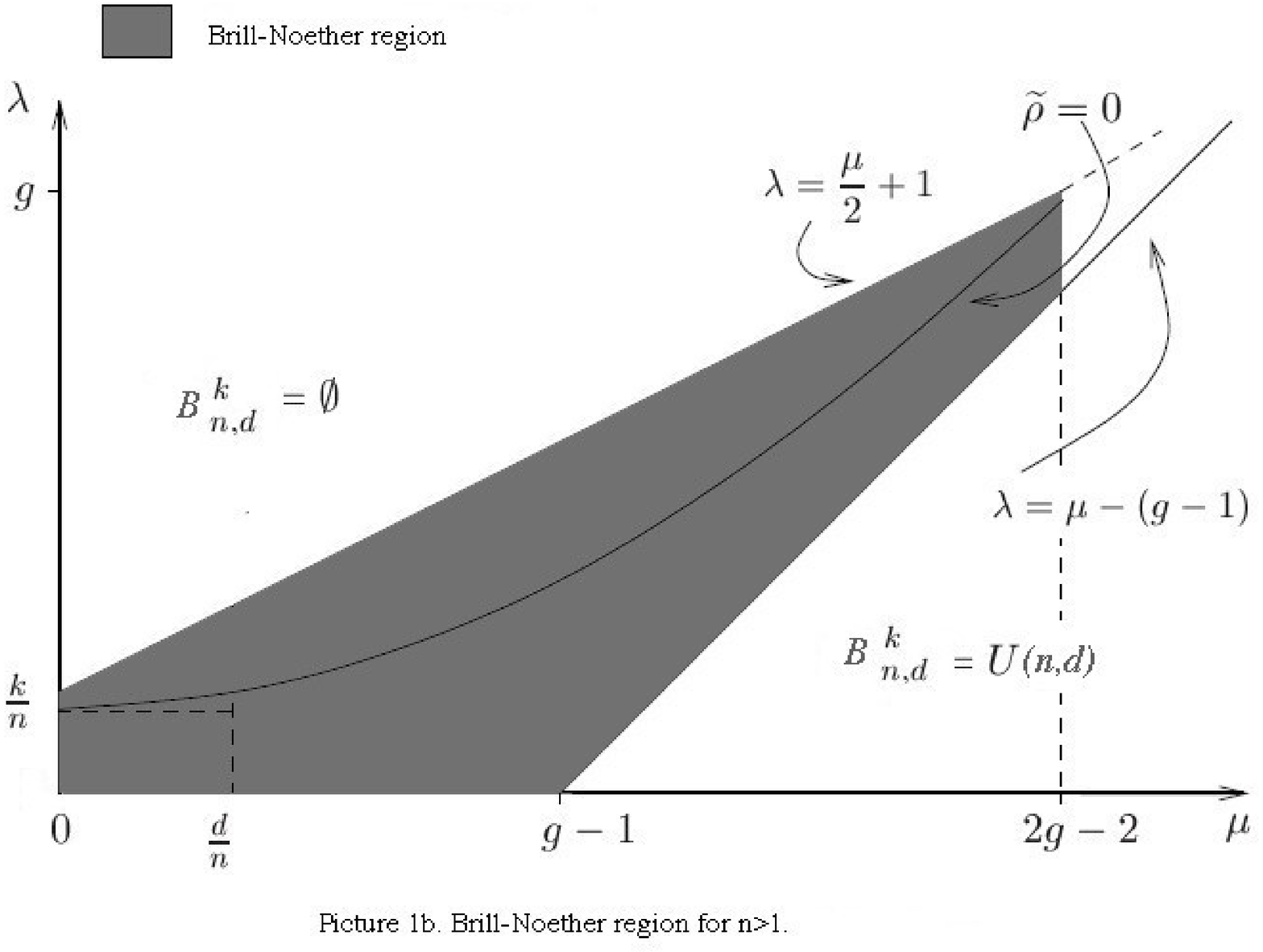,width=#1pt}}

\newcommand{\cpic}[1]{\epsfig{figure=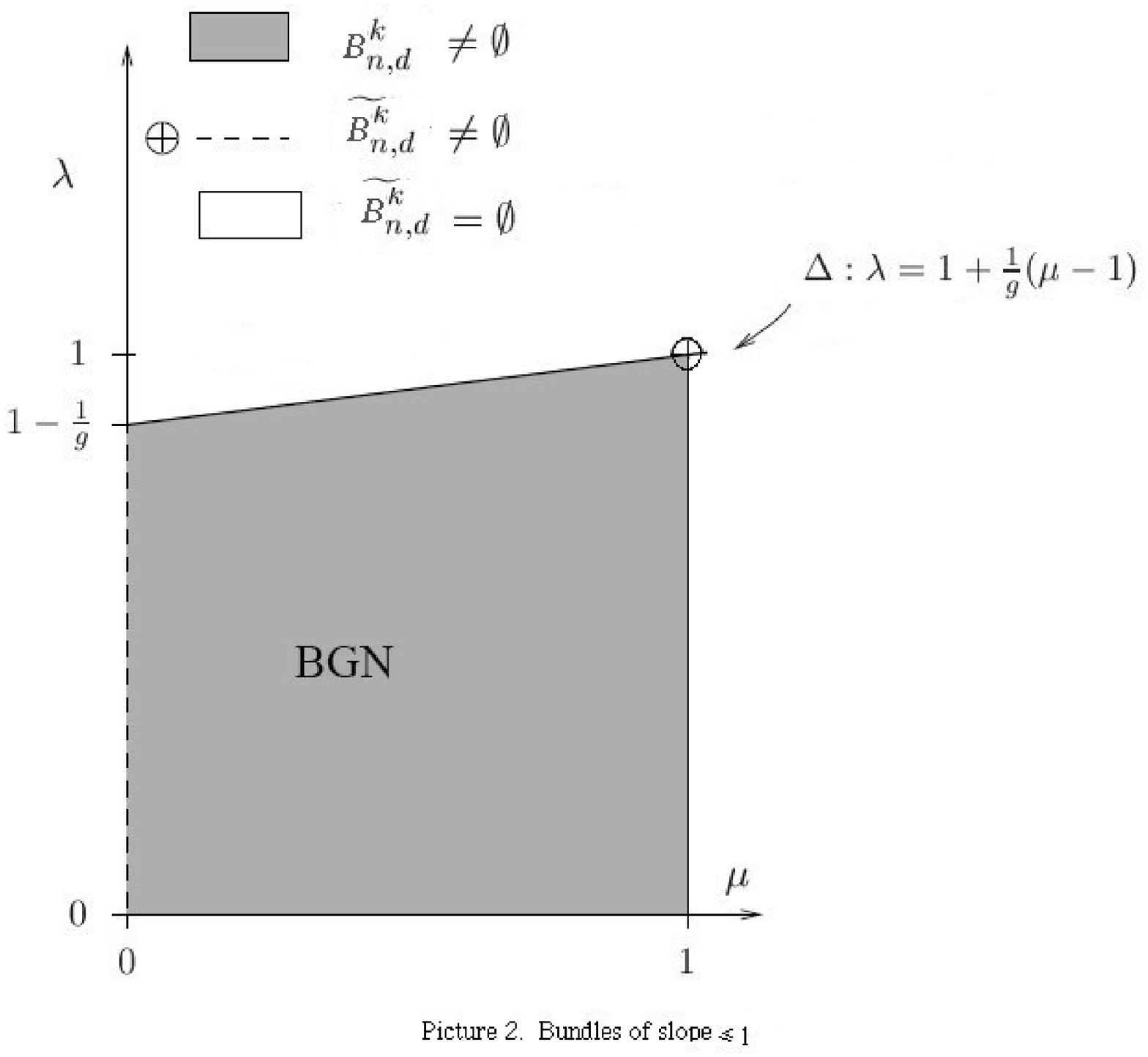,width=#1pt}}

\newcommand{\dpic}[1]{\epsfig{figure=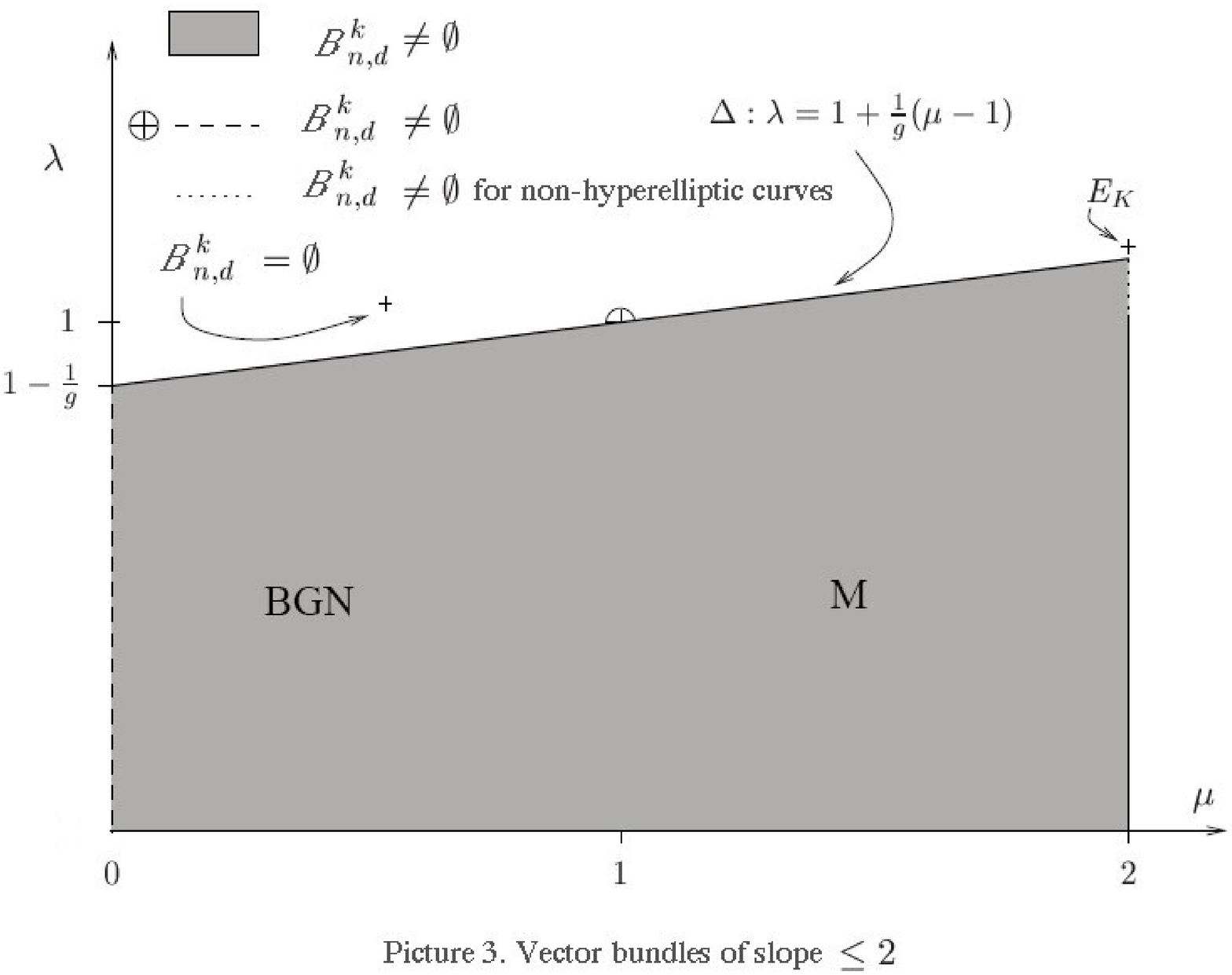,width=#1pt}}

\newcommand{\epic}[1]{\epsfig{figure=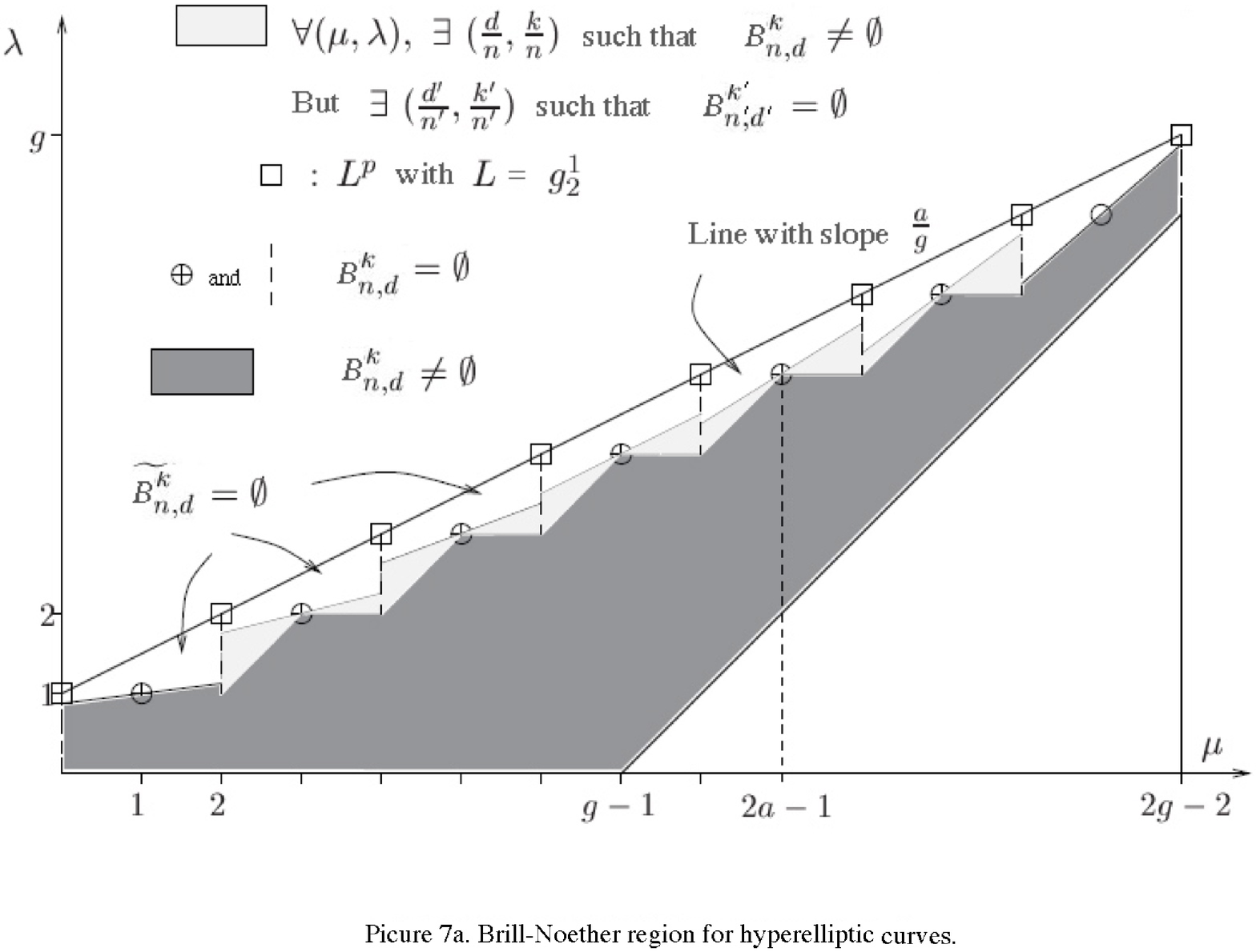,width=#1pt}}

\newcommand{\fpic}[1]{\epsfig{figure=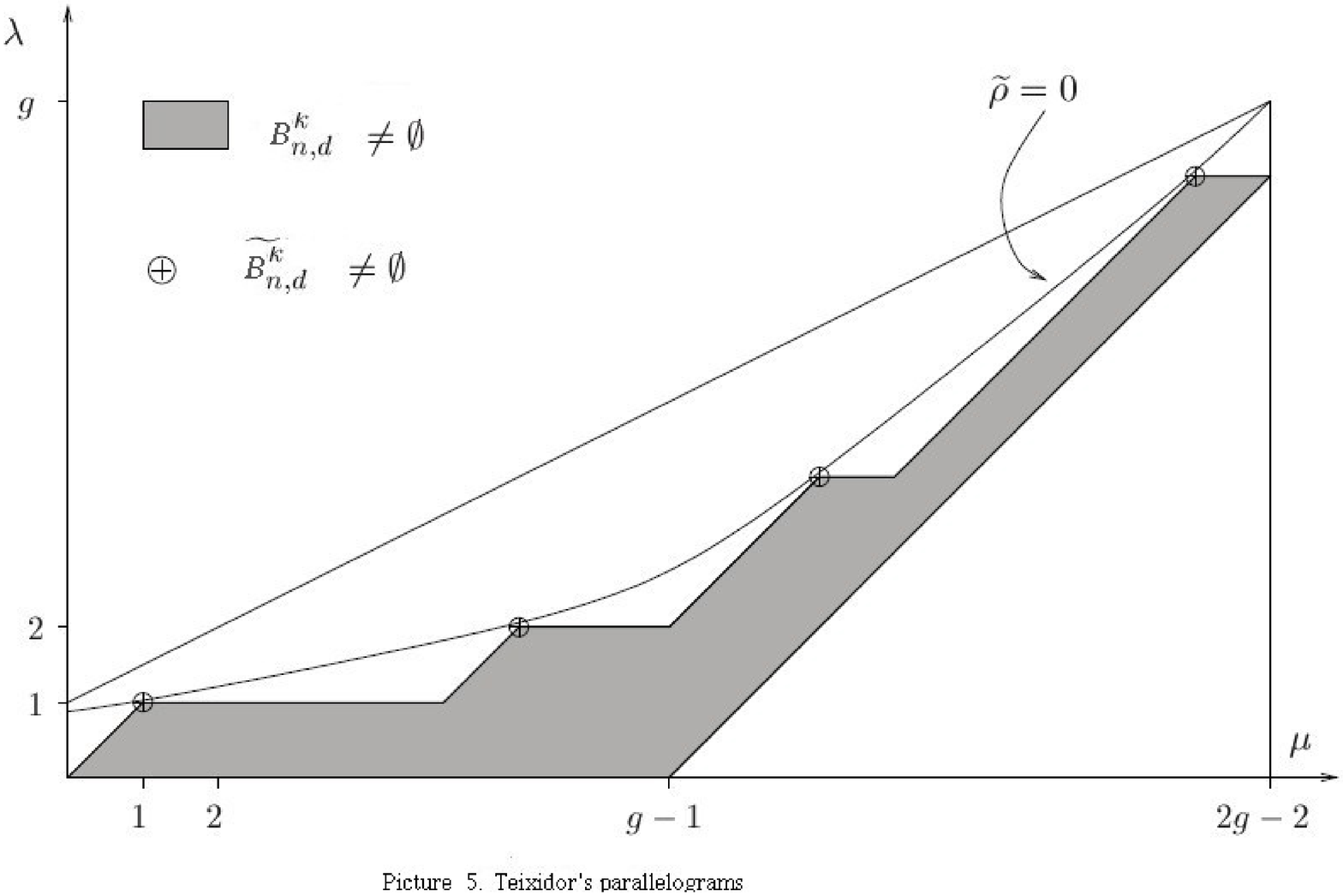,width=#1pt}}

\newcommand{\gpic}[1]{\epsfig{figure=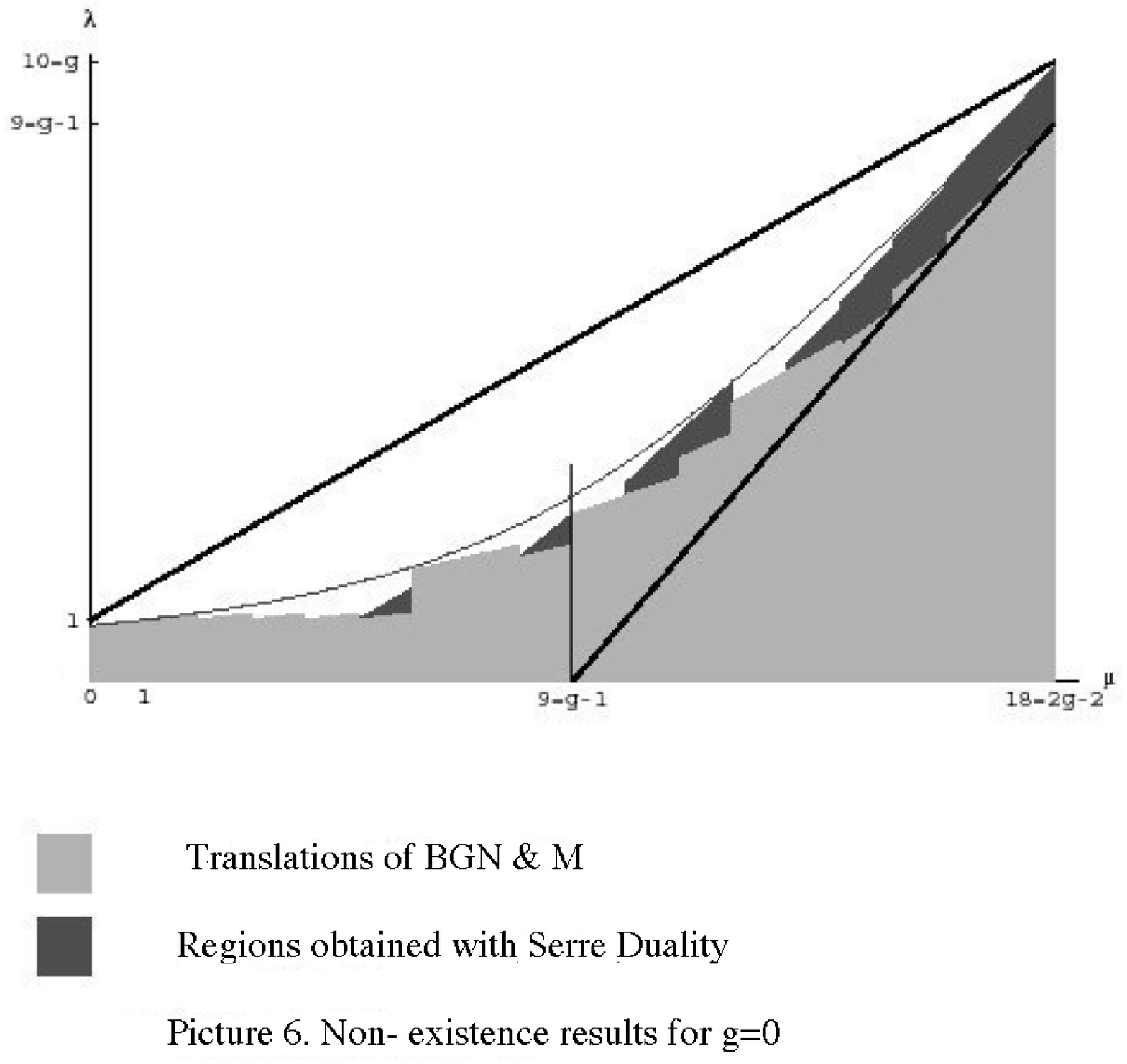,width=#1pt}}

\newcommand{\hpic}[1]{\epsfig{figure=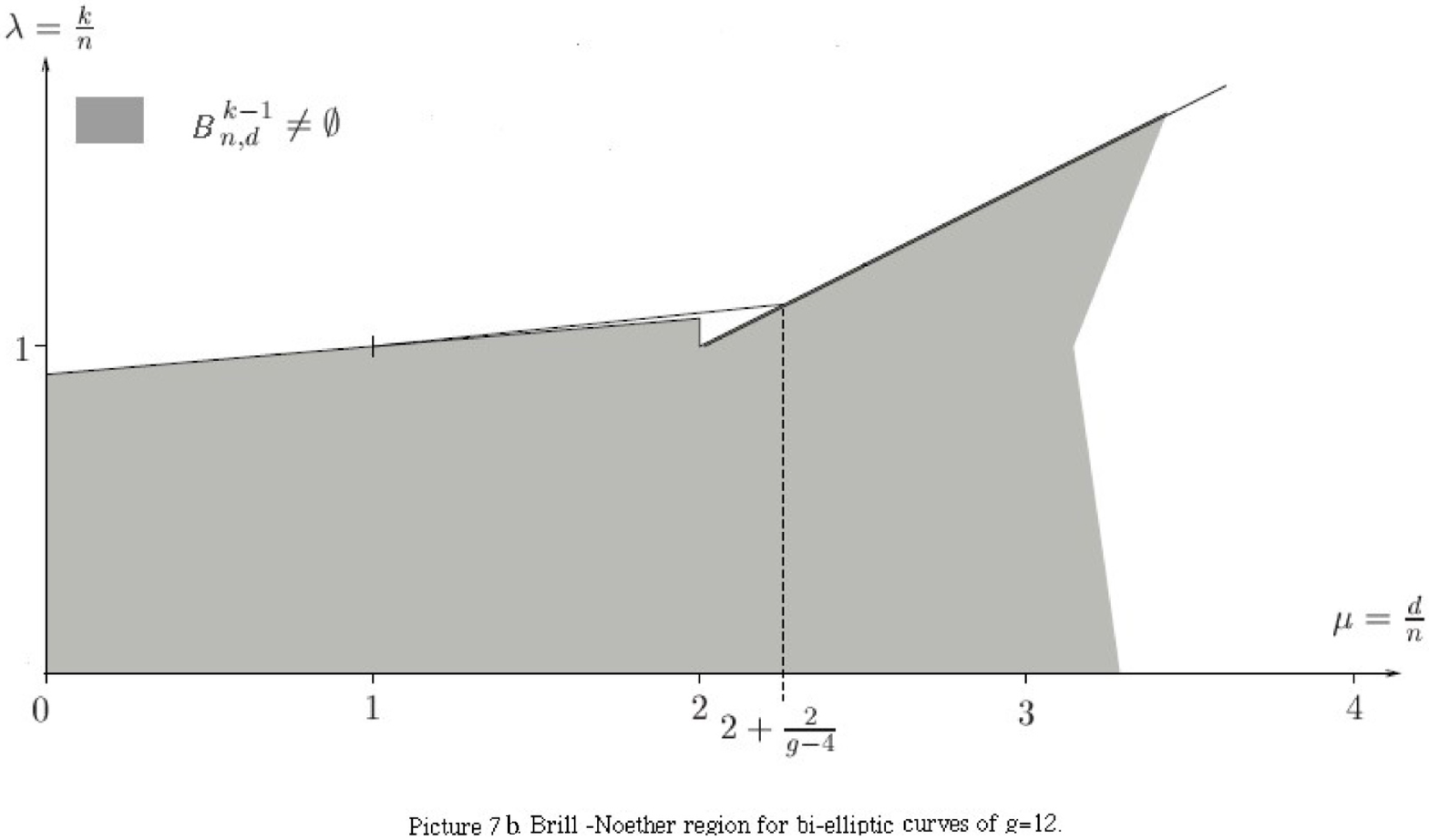,width=#1pt}}

\newcommand{\kpic}[1]{\epsfig{figure=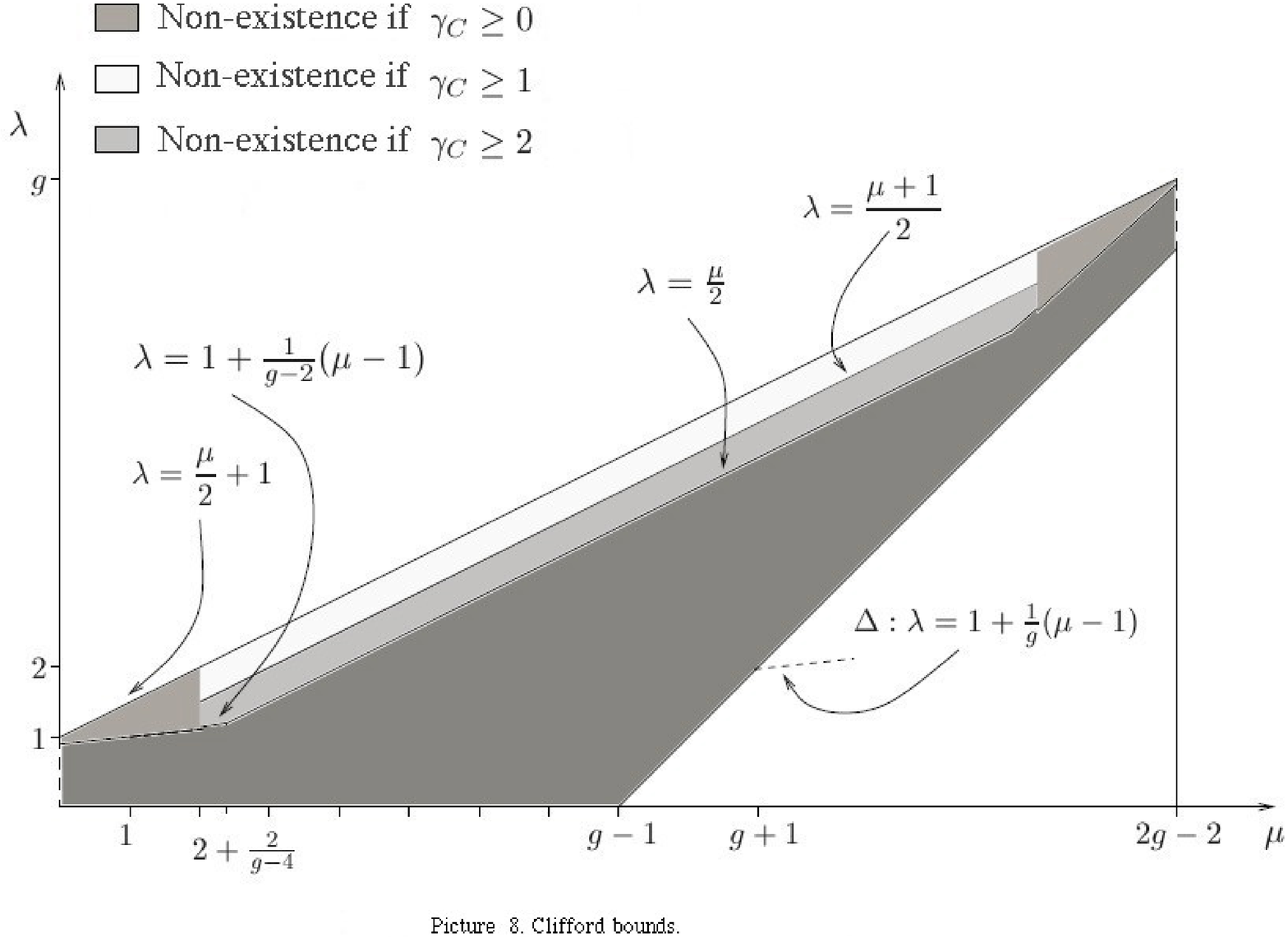,width=#1pt}}

\title{Brill-Noether Theory for  stable vector bundles.}

\author{Ivona Grzegorczyk ; montserrat teixidor i bigas}

\address{Department of Mathematics, California State University Channel Islands, Camarillo, CA 91320}

\email{ ivona.grze@csuci.edu;}

\address{Mathematics Department, Tufts University, Medford MA 02155}

\email{ montserrat.teixidoribigas@@tufts.edu}

\begin{document}

 \begin{abstract} This paper gives an overview of the main results  of Brill-Noether Theory for vector bundles on algebraic curves.

\end{abstract}

\maketitle

\begin{section}{Introduction}

\medskip    Let  $C$  be a projective, algebraic curve of genus  $g$
non-singular most of the time
 and let $ {U(n,d)} $ be the moduli
space  of stable bundles on $ C$  of rank $ n$  and degree $ d$.   A
Brill-Noether subvariety   $B_{n,d}^k$  of $ {U(n,d)} $ is a subset
of  $ {U(n,d)} $  whose points correspond to bundles having at least $k$  independent sections (often denoted by
$W_{n,d}^{k-1}$).
  Brill-Noether theory  describes  the geometry of  $B_{n,d}^k$ and in this paper we present
  the foundations of this theory and an
 overview of  current  results. Even
the basic questions such as when is $B_{n,d}^k$ non-empty, what is its dimension  and whether it is irreducible are subtle and of great interest.  While the topology of $ {U(n,d)} $ depends only on $g,$ the answers to  the above questions turn out to depend on
the algebraic structure of  $C$.   In various cases we have
 definite answers which are valid for a {\it general} curve, but  fail for certain {\it special} curves. This allows to define  subvarieties of the moduli space of curves and contribute to the  study of its geometry. For instance, the first counterexample to  the Harris-Morrison slope conjecture can be defined in such a way  (\cite{FP}).

 Classical Brill-Noether Theory deals with line bundles ($n=1$).
  The natural generalization of most of the results known in this case turn out to be false for some $n\ge 2$ and particular values of $d,k$: the Brill-Noether locus can be empty when the expected dimension is positive and non-empty when it is negative even on the generic curve, it can be reducible on the generic curve with components of different dimensions and the singular locus may be larger than expected.

Moreover, there are very few values of $d,n,k,\ n\ge 2$ for which there is a complete picture of the situation. Nevertheless, for a large set of values of $n,d,k$, there are partial (positive) results.

We will start with a scheme theoretic description of $B_{n,d}^k$ and the study of its tangent space. This justifies its expected dimension and singular locus. We will describe the complete solution for small slopes ($\mu \leq 2 $) (see \cite{BGN}, \cite{Me1}). We will present the results for higher slope from \cite{duke} and \cite{BMNO} and include some stronger results for special curves.
In the last section, we
concentrate on the case of rank two and canonical determinant.

\end{section}

\begin{section}{scheme structure on $B_{n,d}^k$ } We start by giving a scheme structure to $B_{n,d}^k$ . Assume first for simplicity that the greatest common divisor of $n$ and $d$ is one. There exists then a Poincare bundle ${\mathcal E}$ on $C\times U(n,d)$. Denote by $p_1:C\times U(n,d)\rightarrow C$, $p_2:C\times U(n,d)\rightarrow U(n,d)$ the two projections. Let $D$ be a divisor of large degree on $C$, $deg(D)\ge 2g-1-{\frac{d}{n}}$. Consider the exact sequence on $C$ $$0\rightarrow {\mathcal O}_C\rightarrow {\mathcal O}(D)\rightarrow {\mathcal O}_D(D)\rightarrow 0.$$

Taking the pull back of this sequence to the product $C\times U(n,d)$, tensoring with ${\mathcal E}$ and pushing down to $U(n,d)$, one obtains $$0\rightarrow p_{2*}{\mathcal E}\rightarrow p_{2*}{\mathcal E}\otimes p_1^*({\mathcal O}(D))\rightarrow p_{2*}({\mathcal E}\otimes p_1^*({\mathcal O}_D(D)))\rightarrow ...$$

Note now that $p_{2*}({\mathcal E}\otimes p_1^*({\mathcal O}_D(D)))$ is a vector bundle of rank $n$deg$D$ whose fiber over the
point $E\in U(n,d)$ is $H^0(C, E\otimes {\mathcal O}_D(D))$. As $E$ is stable and $K\otimes E^*(-D)$ has negative degree,
$h^0(K\otimes E^*(-D))=0$.
Then, $h^0(C, E(D))=d+n$deg$D+n(1-g)=\alpha$ is independent of the element $E\in U(n,d)$. Hence $p_{2*}({\mathcal E}\otimes
p_1^*({\mathcal O}_D(D)))$ is a vector bundle of rank $\alpha$.
Define, $B^k_{n,d}$ as the locus where the map $$\varphi :p_{2*}{\mathcal E}\otimes p_1^*({\mathcal O}(D))\rightarrow p_{2*}({\mathcal E}\otimes p_1^*({\mathcal O}_D(D)))$$ has rank at most $\alpha -k$.

From the theory of determinantal varieties, the dimension of $B^k_{n,d}$ at any point is at least $\dim U(n,d)-k(n$deg$
D-(\alpha -k))=\rho$ with expected equality. Moreover, provided $B^{k+1}_{n,d} \neq U(n,d)$, the locus where the rank of the
map of vector bundles goes down by one more (namely $B_{n,d}^{k+1} $) is contained in the singular locus, again with expected
equality.

\begin{Thm}{\bf Definition}  The expected dimension of $B^k_{n,d}$ is called the Brill-Noether number and
will be denoted with $$\rho^k_{n,d}=n^2(g-1)+1-k(k-d+n(g-1)).$$

\end{Thm}

When the rank and the degree are not relatively prime, a Poincare bundle on $C\times U(n,d)$ does not exist. One can however take a suitable cover on which the bundle exists and make the construction there.

\bigskip

We want to study now the tangent space to $B^k_{n,d}$ at a point $E$ (see \cite{W}, \cite {BR}).
We first show how the tangent space to $U(n,d)$ at $E$ is identified to $H^1(C, E^*\otimes E)$. Consider a tangent vector to
$U(n,d)$ at $E$, namely a map from $Spec {\bf k}[\epsilon]/\epsilon ^2$ to $U(n,d)$ with the closed point going to $E$. This
is equivalent to giving a vector bundle ${\mathcal E}_{\epsilon}$ on $C_{\epsilon}=C\times Spec {\bf k}[\epsilon]/\epsilon ^2$
extending $E$. Then ${\mathcal E}_{\epsilon}$ fits in an exact sequence $$0\rightarrow E\rightarrow {\mathcal
E}_{\epsilon}\rightarrow E\rightarrow 0.$$ The class of this extension gives the corresponding element in the cohomology
$H^1(C, E^*\otimes E)$.

The bundle ${\mathcal E}_{\epsilon}$ can be described explicitly:
Let $\varphi \in H^1(C, E^*\otimes E)$.
Choose a suitable cover $U_i$ of $C$, $U_{ij}=U_i\cap U_j$.
Represent $\varphi$ by a coboundary $(\varphi _{ij})$ with $\varphi _{ij} \in H^0(U_{ij},Hom(E,E))$.
Consider the trivial extension of $E$ to $U_i\times Spec {\bf k}(\epsilon )/\epsilon ^2$, namely $E_{U_i}\oplus \epsilon
E_{U_i}$. Take gluings on $U_{ij}$ given by $$\begin{pmatrix} Id&0\\  \varphi _{ij}&Id \end{pmatrix}.$$ This gives the vector
bundle ${\mathcal E}_{\epsilon}$.

 Assume now that a section $s$ of $E$
can be extended to a section $s_{\epsilon}$ of the deformation.
Hence there exist  local sections $s'_i \in H^0(U_i, E_{|U_i})$ such that $(s_{|U_i},s_i)$ define a section of ${\mathcal E}_{\epsilon}$.

By construction of ${\mathcal E}_{\epsilon}$ this means that $$\begin{pmatrix} Id&0\\  \varphi_{ij}&Id \end{pmatrix} \begin{pmatrix} s_{|U_i}\\ s'_i \end{pmatrix}= \begin{pmatrix} s_{|U_j}\\ s'_j \end{pmatrix}.$$

The first condition $(s_{|U_i})_{|U_{ij}}=(s_{|U_j})_{|U_{ij}}$ is automatically satisfied as $s$ is a global section. The second condition can be written as $\varphi_{ij}(s)=s'_j-s'_i$.
Equivalently, $\varphi_{ij}(s)$ is a cocycle, namely $$\begin{matrix} \varphi _{ij}\in Ker &( H^1(C, E^*\otimes E)& \rightarrow &H^1(C,E))\\  &\nu _{ij}&\rightarrow &\nu_{ij}(s) \end{matrix}.$$ This result can be reformulated as follows: the set of infinitessimal deformations of the vector bundle $E$ that  have sections deforming  a certain subspace $V\subset H^0(C,E)$  consists of the orthogonal to the image  of the following map

\begin{Thm} {\bf Definition}
The natural cup-product map of sections

$$P_V: V\otimes H^0(C, K\otimes E^*)\rightarrow H^0(C, K\otimes E\otimes E^*).
$$
is called the Petri map. When $V=H^0(C,E)$, we write $P$ for $P_{H^0(C,E)}$.

\end{Thm}

 Note now that if
$E$ is stable, $h^0(C, E^*\otimes E)=1$ as the only automorphisms of $E$ are multiples of the identity.
Hence, $$h^1(C,E^*\otimes E)=h^0(K\otimes E\otimes E^*)=n^2(g-1)+1=\dim U(n,d).$$ Take $V=H^0(C,E)$ above and assume that it
has dimension $k$. By Riemann-Roch Theorem, $h^0(C, K \otimes E^*)=k-d+n(g-1).$
Hence,

$$\rho ^k_{n,d}=h^1(C,E^*\otimes E)-h^0(C,E)h^1(C,E)$$ Equivalently, $B^k_{n,d}$ is non-singular of dimension $\rho ^k_{n,d}$ at $E$ if and only if the Petri map above  is injective.
This result allows us to present the first example of points in the singular locus of $B^k_{n,d}$ that are not in $B^{k+1}_{n,d}$ (see \cite{duke}).

\begin{Thm} {\bf Example.} Let $C$ be a generic curve of genus at least six.
 Let $L$ be a line bundle of degree $g-2$ with at least two independent sections and $E$ a generic extension $$0\rightarrow L\rightarrow E\rightarrow K\otimes L^{-1}\rightarrow 0.$$ Then, $E$ has two sections and the Petri map is not injective.

\end{Thm}

Note first that from the condition on the genus and Brill-Noether Theory for rank one, such an $L$ with two sections exists.

It is easy to check with a moduli count that for a generic extension the vector bundle is stable (see \cite{RT}) and has precisely two sections. Dualizing the sequence above and tensoring with $K$, one obtains an analogous sequence

$$0\rightarrow L\rightarrow K\otimes E^*\rightarrow K\otimes L^{-1}\rightarrow 0.$$ Hence $H^0(C,L)$ is contained both in $H^0(C,E)$ and in $H^0(C, K\otimes E^*)$. Let now $s,t$ be two independent sections of $L$.
Then, $s\otimes t-t\otimes s$ is in the kernel of the Petri map.

For this particular case of rank two and determinant canonical one can check that these are essentially the only cases in which the Petri map is not injective (see \cite{singular}) but not much is known in general.

\end{section}

\maketitle

\begin{section}{Classical Results and Brill-Noether geography}

We start by looking at low genus and rank.

\medskip

{\bf $g=0$.} Grothedieck proved that all bundles on a projective line are decomposable into direct sums of line bundles (see \cite{OSS}), i.e. $E=L_1 \oplus L_2 \oplus...\oplus L_n.$ and there are no stable vector bundles for $n>1$. If $d=nd_1+d_2$, the least unstable vector bundle of rank $n$ and degree $d$ is $${\mathcal O}(d_1+1)^{d_2}\oplus {\mathcal O}(d_1)^{n-d_2}.$$ This vector bundle has an $n^2$-dimensional space of endomorphisms.

The space of pairs $E,V$ where $V$ is a subspace of dimension $k$ of $H^0(C,E)$ is parameterized by the grassmannian of subspaces of dimension $k$ of $H^0(C,E)$ which has dimension $\rho+n^2-1$.

\medskip

{\bf $g=1$. } Atiyah (see [A]) classified the bundles over elliptic

curves: if $(n,d)=1$, the moduli space of stable vector bundles is isomorphic to the curve itself.
If $(n,d)=h>1$, there are no stable vector bundles of this rank and degree and a generic semistable bundle is a direct sum of
$h$ stable vector bundles all of the same rank and degree. From Riemann-Roch, if $d>0$, $h^0(C, E)=d$. Then, the dimension of
the set of pairs $(E,V)$ with $V$ a $k$-dimensional space of sections of $H^0(C, E)$ is the expected dimension $\rho ^k_{n,d}$
of the Brill-Noether locus.

\medskip

{\bf $n=1$.} Line bundles of degree $d$ are classified by the Jacobian variety $J^d(C)=Pic^d(C)$, $ \dim J^d(C)=g.$ All line bundles are stable and $B_{1,d}^{k}$ are subvarieties of $J^d(C)$ and their expected dimension equals $\rho _{1,d}^{k}=\rho = g-k(k-d+g-1)$. To sketch the region  of interest for Brill-Noether theory in the ($d, k$) plane we consider the bounds given by the following Classical results (\cite{ACGH}):

\begin{Thm}{\bf Riemann-Roch theorem for line bundles.} Let $L$ be a line bundle on $C$. Then  $h^0(C,L)-h^1(C,L)=d-(g-1).$

\end{Thm}

\begin{Thm}{\bf Clifford theorem.}  Let $D$ be an effective divisor of degree $d$ and $L={\mathcal {O}}(D)$, $d\leq 2g$.
 Let $h^0(C, {\mathcal O}(D))= k$. Then

$$k< \frac{1}{2}  \deg D+1= \frac{1}{2} \deg L+1.$$ \end{Thm}

\noindent Note, that  when  $\rho =0$ we can solve for $d$

$$d=(k-1)(1+ \frac{g}{k})$$

\noindent i.e. $\rho =0$ is a hyperbola in the $d,k$-plane passing through
(0,1) (corresponding to $O_C$) and through ($2g-2$, $g$) (corresponding to $K_C$).
Points below $\rho = 0$ in the region between the Clifford and the Riemann-Roch lines correspond to $B_{1,d}^k$ of expected
positive dimensions. Serre duality  for line bundles on a curve C given by $h^0(C, {\mathcal O} (D))-h^0(C, {\mathcal
O}(K-D))=d+g+1$ makes two subregions dual to each other, see Picture 1a.

\bigskip

\[\apic{350}\]

The following results give answers to the basic questions of Brill-Noether theory in this case (see [ACGH]).

\begin{Thm}{\bf Connectedness Theorem}. Let $C$  be a smooth curve of genus $g$. Let $d\geq 1, \  k\geq 0$. Assume $\rho = g-k(g-d+k-1) \geq 1$. Then $B_{1,d}^{k}$  is connected on any curve and irreducible on the generic curve.\end{Thm} \begin{Thm}{\bf Dimension Theorem}. Let $C$  be a general curve of genus $g$. Let $d\geq 1, \  k\geq 0$. Assume $\rho = g-k(g-d+k-1) \geq 0$. Then $B_{1,d}^{k}$  is non-empty reduced and of pure dimension $\rho$. If $\rho < 0$,  $B_{1,d}^{k}$ is empty.\end{Thm}

\begin{Thm}{\bf Smoothness Theorem}. Let $C$  be a general curve of genus $g$. Let $d\geq 1, \  k\geq 0$. Assume $\rho = g-k(g-d+k-1) \geq 0$. Then $B_{1,d}^{k}$  is non-singular  outside of $B_{1,d}^{k+1}$ .\end{Thm}

\bigskip

In summary, for rank one the answers to the basic questions are affirmative:
if $\rho \ge 0$, the Brill-Noether
locus is non-empty on any curve and it is connected as soon as $\rho>0$.
It may be reducible, but each component has dimension $\geq\rho.$

Moreover, for the generic curve, the dimension of  the Brill-Noether locus is exactly $\rho$ (with the locus being empty for $\rho <0$), it is irreducible and of dim $\rho$ and Sing $B_{1,d}^{k}$= $B_{1,d}^{k+1}$.

\end{section}

\begin{section}{Problem in higher ranks}

We want to study the non-emptiness of $B_{n,d}^k$. We will make use of two parameters the slope $\mu =\frac{d}{n}$ and a sort of "dimensional slope" $\lambda =\frac{k}{n}$.

To define the area of interest for non-emptiness of $B_{n,d}^k$ in the $\mu, \lambda $ plane, we consider Riemann-Roch Theorem  and Generalized Clifford Theorem  and write them in two new parameters

\begin{Thm}{\bf Riemann-Roch Theorem } The Riemann-Roch line $k  +n(g-1)  \leq d $ { \rm gives} $\lambda +(g-1) \leq \mu$ .

\end{Thm}

\begin{Thm}\label{Clifford} {\bf Clifford's Theorem} (see \cite{BGN}) If $E$ is semi-stable $2k-2n \leq d$ { \rm gives}  $2\lambda -2 \leq \mu .$ \end{Thm}

Using  the expected dimension equation $\rho_{n,d}^{k} = n^2(g-1)+1-k(k-d+n(g-1))$ we define
a new curve $\widetilde {\rho_{n,d}^{k}} = \frac{\rho-1}{n^2} = 0$.
 $\widetilde {\rho_{n,d}^{k}}$ easily translates to
 $\mu $ and $\lambda $ parametization   as

 $$\widetilde {\rho_{n,d}^{k}}=g-1-\lambda (\lambda -\mu +g-1).$$

\medskip

 The Brill-Noether region for $rk \geq 2$ is bounded by the $\widetilde {\rho_{n,d}^{k}}\geq 0$ and Riemann Roch and
Clifford lines. Serre duality describes a ``symmetry" line between the subregions on the left and right, see Picture 1b.

\[\bpic{350}\]

 For $\mu (E) = \frac{d}{n} \leq 1$ the non-emptiness problem of Brill-Noether loci was solved completely by Brambila-Paz,
Grzegorczyk, Newstead (see Picture 2 and [BGN] ) and their results and methods were extended to  $\mu (E) = \frac{d}{n} \leq 2$
by Mercat (see [M1], [M3]).

\[\cpic{350}\]

The main results for $\mu < 2$, $C$ non-singular are stated below.

\begin{Thm}{\bf Theorem }:  Assume $\mu <2$. Then, $B_{n,d}^{k}$ is non-empty iff $d>0$, $n \leq d+(n-k)g$ and $(n,d,k) \neq (n,n,n)$.  If $B_{n,d}^{k}$ is not empty, then it is irreducible, of dimension $\rho_{n,d}^{k}$ and  Sing $B_{n,d}^{k}$= $B_{n,d}^{k+1}.$\end{Thm}

\begin{Thm}{\bf Theorem}: Denote by   $\widetilde {B_{n,d}^{k}}$ the  subscheme of
the moduli space of S-equivalence classes of semistable bundles  whose points correspond to bundles having at least
$k$  independent sections. Then, $\widetilde {B_{n,d}^{k}}$ is non-empty iff either $d=0$ and $k \leq n$ or $d>0$ and $n \leq d+(n-k)g$. If $\widetilde {B_{n,d}^{k}}$ is not empty, then it is irreducible.\end{Thm}

Note: $n=d+(n-k)g$ corresponds to $1=\mu + (1-\lambda )g$ in $(\mu, \lambda)$ coordinates, i.e.
$\lambda =1+ \frac{1}{g}(\mu -1)$ - it is a line which we will denote as $\Delta$. The point $(1,1)$ is on the line and the
$\Delta$ - line is tangent to $\widetilde \rho = 0$ at the point (1,1), see Picture 3.

\bigskip

\[\dpic{350}\]

The main point of the proof in \cite{BGN} is to show that all $E \in B_{n,d}^{k},\ \mu (E) \leq 1$, can be presented as an extension

$$0 \rightarrow O^k \rightarrow E \rightarrow F \rightarrow 0$$ Such an extension corresponds to an element $e=(e_1,...,e_k)\in H^1(F^*)^{\oplus k}$. The condition $n\le d+(n-k)g$ is equivalent to $k\le h^1(F^*)$. If this condition fails, the elements $e_1,...e_k$ are linearly dependent in $H^1(F^*)$ and the extension is partially split. Therefore, the corresponding $E$ cannot be stable.

 For $1<\mu < 2$ consider  the vector bundle defined as  the kernel of the evaluation map of sections of the canonical map
$$(*)\  \  0\rightarrow M_K\rightarrow H^0(K)\otimes{\mathcal O_C}\rightarrow K\rightarrow 0,$$ where $K$ is a canonical
divisor generated by its global sections. As $E\otimes K$ is a stable vector bundle of slope larger than $2(g-1)$,
$H^1(C, E\otimes K)=0$. As $\mu(E)< \mu(M_K^*)=2$, there are no non-trivial maps from $M_K^*$ to $E$, hence $h^0(C, M_K^*\otimes
E)=0$. Then, using the long exact sequence in cohomology of $(*)$ tensored with $E$, one can see that the inequality $n \leq
d+(n-k)g$ still holds in this case.

These results were extended to  $\mu = 2$ as follows:

\begin{Thm}{\bf Theorem } : If $C$ is nonhyperelliptic then $B_{n,d}^{k}$ is non-empty iff  $n \leq d+(n-k)g$  or
$E\simeq M_K^*$.

If $C$ is hyperelliptic, then $B_{n,d}^{k}$ is non-empty iff  $k \leq n$ or $E\simeq g^1_k$.
There is at least one component of the right dimension equal to $\rho_{n,d}^{k}$.

\end{Thm}

\begin{Thm}{\bf Comments.}\end{Thm}

1. The material in this section provides many examples with  $\rho_{n,d}^{k} > 0 $ and $B_{n,d}^{k}= \phi$.
 These are unexpected results.

2. Note that $\widetilde{B_{n,0}^{k}}$ is not a closure of $B_{n,0}^{k}=\phi$, i.e semistable bundles with sections exist when
stable ones do not.

3. Recent results on coherent systems show that when $1\leq \mu \leq2$ there are no other components of $ B_{n,d}^{k}$, (see [BGMMN]).

4. For $\mu (E) \leq 2$, the $\Delta$- line that is tangent to $\widetilde{\rho_{n,d}}^{k}=0 $ at the point $(1,1)$ is the nonemptiness boundary (even though it is completely in the Brill-Noether region).  Does this indicate existence of similar boundaries for other $\mu $'s?

\end{section}

\maketitle

\begin{section}{Non-emptiness of Brill-Noether loci for large number of sections}

In this section, we expose the main known results about existence of vector bundles with sections mostly when $k>n$, as in the case $k\le n$ better results can be obtained with the methods of the previous section.

\begin{Thm}\label{existence} {\bf Theorem} Let $C$ be a generic curve of genus $g\ge 2$. Let $d,n,k$ be positive integers with $k>n$. Write $$d=nd_1+d_2,\ k=nk_1+k_2,\ d_2<n, k_2<n$$ Then, $B^k_{n,d}$ is non-empty and has one component of the expected dimension $\rho$ if one of the conditions below is satisfied:

$$g-(k_1+1)(g-d_1+k_1-1)\ge 1,\ 0\not= d_2\ge k_2$$ $$g-k_1(g-d_1+k_1-1)> 1,\  d_2= k_2=0$$ $$g-(k_1+1)(g-d_1+k_1)\ge 1,\  d_2< k_2 .$$ \end{Thm}
See Picture 5 for graphic presentation of this result in the Brill-Noether region.

\[\fpic{350}\]

We sketch here two methods of proof. The first one is based in degeneration techniques. Assume given a family of curves $${\mathcal C}\rightarrow T$$ and construct the Brill-Noether locus for the family in a way similar to what was done in section 2 for a single curve $${\mathcal B}^k_{n,d}=\{ (t,E) |t\in T,\ E\in B^k_{n,d}({\mathcal C}_t) \}.$$

The dimension of this scheme at any point is at least $\dim T+\rho ^k_{n,d}$. Assume that we can find a point $(t_0,E)$ such that $\dim( B^k_{n,d}({\mathcal C}_{t_0} ))_E=\rho^k_{n,d}$. Then $\dim {\mathcal B}^k_{n,d}\le \rho ^k_{n,d} + \dim T$ and therefore we have equality.

The dimension of the generic fiber of the projection ${\mathcal B}^k_{n,d}\rightarrow T$ is at most the dimension of the fiber over $t_0$, namely $\rho ^k_{n,d}$. But it cannot be any smaller, as the fiber over a point $t\in T$ is $B^k_{n,d}({\mathcal C}_t)$ which has dimension at least $\rho ^k_{n,d}$.  Hence there is equality which is the result we are looking for.

The main difficulty of course is finding a particular curve and vector bundle for which the result works. The curves that seem to work best are obtained by taking $g$ elliptic curves $C_1,...,C_g$ each glued to the next by a chain of rational components. There are three technical problems that need to be solved with this approach:

1) The limit of a vector bundle on a non-singular curve is a torsion-free, not necessarily locally-free sheaf on a nodal curve.

2) The definition of stability cannot be applied directly to the reducible situation to produce moduli spaces of sheaves on reducible curves.

3) It is not clear what one should mean by the limit of the $k$-dimensional space of sections.

The first problem can be solved by blowing up some of the nodes and replacing the curve by a curve with some more components (see \cite{G}, \cite{EH}).

The second problem was solved in two different ways by Gieseker and Seshadri(see \cite{G}, \cite{S}). Gieseker considered the
concept of Hilbert stability: a vector bundle and a space of sections that generate it give rise to a map from the curve to the
Grassmannian.

He then took the Hilbert scheme of non-singular curves in the Grassmannian and considered its closure. Seshadri on the other hand, generalized the concept of stability to torsion-free sheaves on reducible curves: given a curve $C$ with components $C_i$, take a positive weight $a_i$ for each component with $\sum a_i=1$. Define the $a_i$-slope of a torsion-free sheaf $E$ as $\frac{\chi (E)}{\sum a_irank(E_{|C_i})}$. As usual, a sheaf is stable if its slope is larger than the slope of any of its subbundles. Seshadri then constructed a moduli space for $a_i$-stable sheaves.

 The third problem is probably the trickiest and requires a generalisation of the techniques of Eisenbud and Harris for
limit linear series (\cite{EH}) to higher rank (see \cite{duke}).

Once these three obstacles have been overcome, one can use the knowledge of vector bundles on elliptic curves to obtain the results. The main advantage of this method is that it can still be used for numerical conditions different from those in \ref{existence} to improve the results in particular situations (see next section).

\bigskip

A second approach was taken by Mercat (\cite{Me2}). He proved the following

\begin{Thm} {\bf Lemma}. Let $L_1,...L_n$ be line bundles of the same degree $d'$ not isomorphic to each other and $F$ either a vector bundle of slope higher than $d'$ or a torsion sheaf. Then, there exists a stable $E$ that can be written as an extension $$0\rightarrow L_1\oplus...\oplus L_n\rightarrow E\rightarrow F\rightarrow 0.$$ \end{Thm}

One can then deduce \ref{existence} by taking line bundles $L_i$ with lots of sections. For example, if $d_2<k_2$, take $L_i\in B^{k_1+1}_{1,d_1}$ (which are known to exist by classical Brill-Noether Theory) and $F$ a torsion sheaf of degree $d_2$.

 This method
does not give the fact that one can get a component of dimension $\rho$. In fact the set of $E$ as described above when the $L_i$ have plenty of sections form a proper subset of any component of $B^k_{n,d}$. On the other hand,  for certain choices of $d,n,k$ one can obtain an $E$ with $k$ sections for which $\rho^k_{n,d}<0$.

Namely, one obtains {\bf non-empty} Brill-Noether loci with {\bf negative} Brill-Noether number for every curve of a given genus.

\end{section}

\begin{section}{Further results}

The results for small slopes can be used to obtain non-emptiness for larger slopes (and small number of sections) as follows:

The "translation method" based on twisting and dualizing vector bundles with known properties used in [BMNO]  leads to extension of the previous results for a general curve $C$ to give the following

results:

\begin{Thm}

 {\bf Theorem} (\cite{BMNO} ) $ {B_{n,d}^{k}} $ is not empty if the following conditions are satisfied $d=nd'+d''$ with
$0<d''<2n$, $1<s<g$, $d' \geq \frac{(s-1)(s+g)}{s}$, $n \leq d''+(n-k)g$,

$(d'',k) \neq (n,n).$

\end{Thm}

The above  gives polygonal regions below the $\widetilde{\rho_{n,d}}^{k}=0 $ curve that sometimes go beyond
Teixidor parrallelograms, see Picture 6. Note though that these are only valid for small values of $k$ ($k\le \frac{d''}{g}
+\frac{n(g-1)}{g}<2n$).

For special values of $k,n,d$, further results can be obtained.

\[\gpic{350}\]

More can be said of course when we restrict the values of the parameters. For example in rank two the best result so far seems to be the following:

\begin{Thm}{\bf Theorem}(see \cite{T5}) If $k\ge 2$ and $k-d+2(g-1)\ge 2$, then, $B^k_{2,d}$ is non-empty and has one component of the right dimension on the generic curve for $\rho \ge 1$ and odd $d$ or $\rho \ge 5$ and $d$ even.

\end{Thm}

 If $g \geq 2$ and $3 \leq d \leq 2g-1$ then $B_{2,d}^{2} \neq \phi$, and the locus is irreducible, reduced of dimension $\rho =2d-3$ and if $d>g+1$  Sing $B_{2,d}^{2}=B_{2,d}^{3} \cup F$  ($\dim F=2g-6$).

If $g \geq 3$ and  $0<d<2g$ then $B_{2,d}^{3} \neq \phi$ iff $\rho_{2,d}^{3} \geq 0.$\\ Morover if $\rho >0$, $B_{2,d}^{3}$ is irreducible, reduced of $dim=\rho$.

It is well known that special curves do not behave as generic curves from the point of view of classical Brill-Noether Theory. The same is of course true for higher rank. Here is a sampling of results:

\[\epic{350}\]

Results or hyperelliptic curves (see Picture 7a above for summary and [BGMNO] for details)  give examples of non-empty loci for
values of
$\lambda,
\mu$ very close to the Clifford line.

Bielliptic curves also behave differently from generic curves and known results are presented on Picture 7b below and stated in
the following

\begin{Thm}{\bf Theorem} (Ballico) Let $C$ be a bielliptic curve (i.e.there exists a two to one map $f: C\rightarrow C'$ where $C'$ is an elliptic curve. Then

If $h^0(E) \leq \frac{d}{2}$ ($\lambda \leq \frac{\mu}{2}$) then $\widetilde{B_{n,d}^{k}} \neq \phi$.  \\

If $h^0(E) < \frac{d}{2}$ ($\lambda < \frac{\mu}{2}$) then $B_{n,d}^k \neq \phi$. \end{Thm}

\bigskip

\[\hpic{350}\]

\bigskip

Denote by $B^k_{n,L}$ the locus of stable vector bundles of rank $n$ and determinant $L$ with at least $k$ sections.
One has the following result (see [V]):

\begin{Thm} {\bf Theorem}. Le $C$ be a curve lying on a $K3$ surface. Then

$B_{2,K}^{k} \neq \phi$ if $k \leq [\frac{g}{2}] +2$.\end{Thm}

The $B_{2,K}^{[\frac{g}{2}]+2}$  locus is known to be empty on the generic curve. Hence curves lying on $K3$
surfaces are exceptional from the point of view of Brill-Noether Theory of rank two and determinant
canonical. This is particularly relevant as generic curves on $K3$'s are well-behaved from the point of
view of classical Brill-Noether Theory. Therefore, higher rank Brill-Noether provides a new tool to study
the geometry of ${\mathcal M}_g$ (see also some very interesting related results in \cite{AN}).

\end{section}

\begin{section}{Generalized Clifford bounds}

Clifford's Theorem for line bundles was generalized to semistable vector bundles by Xiao (\ref{Clifford}). The first extension of this bound was given by Re in the following :

\begin{Thm} {\bf Theorem}. (see \cite{R})  Let $C$ be a nonhyperelliptic curve , $E$ a stable vector bundle on $C$ such that
$1\leq\mu\leq 2g-1$ then $h^0(E) \leq \frac{d+n }{ 2}, $ (or $\lambda \leq \frac{\mu +1}{2}$).\end{Thm}

This inequality moved the Clifford bound lower in the $(\mu , \lambda )$ plane towards the $\widetilde \rho$ curve, giving more non-existence of stable vector bundles results on nonhyperelliptic curves.

To extend this idea one makes the following definition:

\begin{Thm} {\bf Definition of Clifford index}. The Clifford index, denoted by  $\gamma_C$, is the maximal integer such that
$h^0(C, L)\leq \frac{{\deg L - \gamma_C} }{ 2} +1$ for line bundles $L$ such that $h^0(L) \geq 2$ and $h^1(L) \geq 2$.\end{Thm}

\begin{Thm} {\bf Comments}. \end{Thm} 1.  $\gamma_C=0$ iff $C$ hyperelliptic.

2. $\gamma_C=[\frac{g-1 }{ 2}]$ for $C$ of genus $g$ generic.

The following theorem in its generality is still a conjecture:

\begin{Thm} {\bf Generalized Clifford Bound Theorem}.  For $g\geq 4$ and $E$ semistable,
$h^0(E) \leq \frac{d-\gamma_Cn}{ 2}+n$ for $\gamma_Cn+2\leq \mu (E)\leq 2g-4-\gamma_Cn$. \end{Thm}

However, the  non-existence is established in the Theorem 7.5 below, see Picture 8 or graphic interpretation and [M4] for
details.

\begin{Thm} {\bf Theorem.}
 If \
$2
\leq
\mu (E)
\leq 2+\frac{2}{ g-4}$, then \
$h^0(E)\leq n+ \frac{1}{ g-2} (d-n)$.

If\  $2 + \frac{2}{g-4} \leq \mu(E) \leq 2g-4-\frac{2}{ g-4}$, then \
$h^0(C, E) \leq \frac{d}{2}$.

\end{Thm}

\[\kpic{350}\]

\bigskip

For a Petri curve (i.e. a curve for which the Petri map is injective in rank one) of genus 10, the Clifford index is 4. From 7.4 for vector bundles of rank 2 we get $h^0(E) \leq \frac{d-8}{ 2}+2$ for $6 \leq \mu (E) \leq 12. $ For example for $d=12$, the point corresponding to $B_{2,12}^{4}$ represents a "corner point" for the Clifford bound inequalities, (see [GMN]).

\begin{Thm} {\bf Theorem}. Let $C$ be a general curve of genus 10. then $B_{2,d}^{4}\neq \phi $ if and only if $d\geq 13$. Hovever, there exists a curve for which $B_{2,12}^{4}\neq \phi $\end{Thm} This result shows again, that even the basic properties of $B_{n,d}^{k}$ may depend on the curve structure, (see [GMN]).

\end{section}

\medskip

\begin{section}{The case of rank two and canonical determinant}

Assume now that $E$ is a vector bundle of rank two and canonical determinant. In this section, we shall assume that the characteristic of the ground field is different from two. We want to consider the set $$B_{2,K}^k=\{ E\in U(2,K)|h^0(C,E)\ge k \}.$$

\begin{Thm}{\bf Proposition} (see \cite{BF}, \cite{M1}) The set $B_{2,K}^k$ can be given a natural scheme structure. Its expected dimension is $$\dim \ U(2,K)-{k+1\choose 2}.$$

The tangent space to $B_{2,K}^k$ at a point $E$ is naturally identified to the orthogonal to the image of the symmetric Petri map $$S^2(H^0(C,E))\rightarrow H^0(C,S^2(E)).$$ \end{Thm}

\begin{proof}

The proof uses the following fact:

 \begin{Thm}\label{antisimm.} {\bf Lemma}(see \cite{H, Mu}).

 Let $\alpha$ be an antisymmetric non-degenerate form on a vector  bundle ${\mathcal F}$ of rank $2n$ on a variety $X$. Let  $A$  and $A'$ be isotropic subbundles of rank $n$. Then, $$\{ x\in X| \dim(A_x\cap A'_x)\ge k ,\ A_x\cap A'_x\equiv k (mod \ 2)\}$$ has codimension at most ${k+1\choose 2}$ in $X.$ \end{Thm}

 The construction that follows is a direct generalization of the one in \cite{H}.
 As in section 2, let us assume for simplicity that there exists  a Poincare bundle ${\mathcal E}$ on $C\times U(2,K)$.
Denote by $p_1:C\times U(2,K)\rightarrow C$, $p_2:C\times U(2,K)\rightarrow U(2,K)$ the two projections. Let $D$ be a divisor
of large degree on $C$. Define $${\mathcal F}= p_{2*}(({\mathcal E}\otimes p_1^*{\mathcal O}_C(D))/({\mathcal E}\otimes
p_1^*{\mathcal O}_C(-D)).$$ This is a vector bundle on $U(2,K)$ of rank 4deg$D$. Define $N=2$deg$D.$

Define an antisymmetric form $\alpha $ in ${\mathcal F}$ by $\alpha (s_i,t_j)=\sum_{P\in D}Res(s_i\wedge t_j)$ where $s_i\wedge t_j$ is taken as a section of $K(D)$ and $Res$ denotes the residue.

Note now that

$$A= p_{2*}({\mathcal E}/({\mathcal E}\otimes p_1^*{\mathcal O}_C(-D)))$$ is a subbundle of rank $N$ that is
isotropic (as we are considering holomorphic sections only, the residues are zero). Also $$A'= p_{2*}({\mathcal E}\otimes
p_1^*{\mathcal O}_C(D))$$ is an isotropic subbundle: by the residue Theorem, the sum of the residues at all the poles is zero.
The intersection of the fibers of these two subbundles at a point $E$ is $H^0(C,E)$. Then the locus $$B^k_{2,K}=\{ x\in U(2,K)|
\dim (A_x\cap A'_x)\ge k \}.$$

Therefore, by \ref{antisimm.}, the expected dimension of $B_{2,K}^k$ is as claimed.

Let us now turn to the computation of the tangent space. We first need to identify the tangent space to $U(2,K)\subset U(2,2g-2)$ as a subset of $H^1(C,E^*\otimes E)$.
Denote by $Trl(E)$ the sheaf of traceless endomorphisms of $E$. We have a natural decomposition of $E^*\otimes E$ into direct sum

$$\begin{matrix} E^*\otimes E&\rightarrow & {\mathcal O}_C\oplus Trl(E)\\

                   \varphi&  \rightarrow &(\frac{Tr(\varphi)}{ 2}Id,\varphi-\frac{Tr(\varphi)}{ 2}Id)\\

                   \end{matrix}$$

the inverse map being given by the two natural inclusions. We can identify the tangent space to $U(2,K)$ with $H^1(C, Trl(E))$. This can be seen with an infinitesimal computation as in section 2. With the notations in section 2, let $\varphi_{ij}$ be a representative of an element in $H^1(C,E^*\otimes E)$. Take a local basis $e_1,e_2$ for the fiber of $E$ on the open set $U_i$. Then, $e_1\wedge e_2$ is a basis for $\wedge^2 E$ in $U_i$. We want to see which deformations ${\mathcal E}_{\epsilon} $ of $E$ preserve the determinant. Note that $e_1\wedge e_2$ glues in $U_j$ with $$(e_1+\epsilon \varphi _{ij}(e_1))\wedge (e_2 +\epsilon\varphi _{ij}(e_2))=e_1\wedge e_2+\epsilon (\varphi _{ij}(e_1)\wedge e_2 +e_1\wedge \varphi _{ij}(e_2))$$ Writing $\varphi_{i,j}$ as a square matrix in terms of the basis $e_1,e_2$, $$\varphi _{i,j}=\begin{pmatrix}a_{11}&a_{12}\\a_{21}&a_{22}\\\end{pmatrix}$$

one can check that

$$\varphi_{ij}(e_1)\wedge e_2 +e_1\wedge \varphi _{ij}(e_2)= (a_{11}e_1+a_{12}e_2)\wedge e_2+e_1\wedge (a_{21}e_1+a_{22}e_2)=$$$$(a_{11}+a_{22})(e_1\wedge e_2)= Tr(\varphi_{i,j})(e_1\wedge e_2)$$ Therefore, the deformation preserves the determinant if and only if the trace is zero. Hence, the tangent space to $U(2,K)$ at $E$ is $H^1(C,Trl(E))$.

From the computations in section two, the deformations that preserve a certain subspace of sections correspond to the orthogonal to the image of the Petri map. We need to compute the intersection of this orthogonal with $H^1(C,Trl( E))$.

From $\wedge ^2(E)\cong K$, we obtain the identification $E\cong E^*\otimes K\cong Hom(E,K)$ given by $$\begin{matrix} E&\rightarrow & Hom(E,K)& & & \\

                  e&\rightarrow &\psi_e:& E&\rightarrow &K\\

 &             &       & v&\rightarrow & e\wedge v\\

 \end{matrix}.$$

 Similarly, one has the identification

  $$E\otimes E \cong E^*\otimes K\otimes E\cong Hom(E, K\otimes E).$$

  Consider the composition of this identification with the direct sum decomposition of $E^*\otimes E$ tensored with $K$ $$E\otimes E\cong E^*\otimes K\otimes E\cong  K\oplus (K\otimes Trl(E)).$$

 The projection on the first summand is (up to a factor of two)  $$e\otimes v\rightarrow Tr(\psi _e\otimes v).$$  Computing locally, one checks that this is the determinant  $$e\otimes v\rightarrow e\wedge v.$$ Therefore, $H^0(C,K\otimes Trl(E))$ can be identified to the kernel of this map, namely $H^0(C,S^2(E))$.

With these identifications, the composition of the Petri map with the projection onto $H^0(C,S^2(E))$ becomes $$H^0(C,E)\otimes H^0(C,E)\rightarrow H^0(C,E\otimes E)\rightarrow H^0(C,S^2(E)).$$ This map vanishes on $\wedge ^2 H^0(C,E)$ and therefore gives rise to the symmetric Petri map $$S^2(H^0(C,E))\rightarrow H^0(C,S^2(E)).$$ Then, the orthogonal to the image of this map is identified to the tangent space at $E$ to $B^k_{2,K}$.

\end{proof}

As in the case of non-fixed determinant, the map being injective is equivalent to $E$ being a non-singular point on a component of the expected dimension of $B^k_{2,K}$. This is the subject of \cite{Mukai}.

\medskip

In terms of existence, one has the following result

\begin{Thm} {\bf Theorem} (see \cite{canonic} ) Let $C$ be a generic curve of genus at least two. If $k=2k_1$, then $B_{2,K}^k$
is non-empty and has a component of the expected dimension if $g\ge k_1^2, \ {\textrm if}\ k_1>2; g\ge 5 , \ {\textrm if}\
k_1=2; g\ge 3, \ {\textrm if}\ k_1=1.$ If $k=2k_1+1$, $B_{2,K}^k$ is non-empty and has a component of the expected dimension
if $$g\ge k_1^2+k_1+1 . $$ \end{Thm}

This result was proved using degeneration techniques. Very similar results were obtained in \cite{P} computing the cohomology class of the locus inside $U(2,K)$ and showing that it is non-zero.

Note that for small values of $k$, the expected dimension of $B_{2,K}^k$ is larger than the expected dimension of $B_{2,2g-2}^k$.
This allows us to produce some more counterexamples to the expected results. Here we state only the case of an odd number of
sections:

\begin{Thm} {\bf Example} Let $C$ be a generic curve of genus $g\ge 2$.
If $k_1^2+k_1+1\le g< 2k_1^2+k_1$, then $B_{2,2g-2}^{2k_1+1}$ is reducible and has one component of
dimension larger than expected.

\end{Thm}

\end{section}

\end{document}